\documentclass[12pt]{article}

\usepackage[a4paper, total={6in, 9in}]{geometry}

\usepackage{mathptmx}
\usepackage{times}
\usepackage{verbatim}
\usepackage{amsfonts}
\usepackage[normalem]{ulem}

\usepackage{amsmath,amsfonts,amssymb,amsthm}
\usepackage{bm}

\usepackage{hyperref}
\usepackage{color}
\usepackage{booktabs}
\newcommand\gen{{\mathrm {Gen}}}

\usepackage{tikz}
\usepackage{graphicx} 
\usepackage{float}
\tikzset{
    ns/.style={shape=circle, inner sep=1pt, color=black},
    nsf/.style={shape=circle, inner sep=1.2pt, color=black, fill=black} 
}

\usepackage[numbers,sort&compress]{natbib}

\allowdisplaybreaks
\def\qed{\nopagebreak\hfill{\rule{4pt}{7pt}}}
\def\proof{\noindent {\it{Proof.} \hskip 2pt}}
\newtheorem{theorem}{Theorem}[section]

\newtheorem{proposition}[theorem]{Proposition} 
 
\numberwithin{equation}{section}

\begin{document}
 
\begin{center}

 {\Large\bf A Grammatical Calculus for the Ramanujan  Polynomials }
\\[20pt]

William Y.C. Chen$^1$,   Amy M. Fu$^2$ and Elena L. Wang$^3$

\vskip 3mm

$^{1,3}$Center for Applied Mathematics, KL-AAGDM \\
Tianjin University\\
Tianjin 300072, P.R. China

\vskip 3mm

$^{2}$School of Mathematics\\
Shanghai University of Finance and Economics\\
Shanghai 200433, P.R. China

\vskip 3mm

Emails: { $^1$chenyc@tju.edu.cn, $^{2}$fu.mei@mail.shufe.edu.cn,
$^3$ling\_wang2000@tju.edu.cn
}
   
\end{center}

  \begin{center}
  {\it Dedicated to Krishnaswami Alladi for his leadership and friendship}
  \end{center}

\vspace{0.6cm}

\centerline{\bf Abstract}

   The Ramanujan polynomials arise in three
   intertwined contexts. As remarked by Berndt-Evans-Wilson, no combinatorial perspective seems to be
    alluded to in the original definition
   of Ramanujan.
   On a different stage, 
   Dumont-Ramamonjisoa uncovered a combinatorial
   structure underneath an equation also
   considered by Ramanujan. Around the same time, 
   Shor came up with the same
   construction as a
   refinement of the classical formula of 
   Cayley for trees. 
   We present a labeling scheme for 
   rooted trees by employing an extra label marking improper edges.
   Harnessed by this grammar, we develop a 
   grammatical calculus for
   the Ramanujan polynomials heavily
   relying on the constant properties.
   Moreover, we
   provide a grammatical formulation of 
a correspondence that leads to the
recurrence relation  due to 
Berndt-Evans-Wilson and Shor.

\vspace{0.3cm}



\section{Introduction}

Speaking of the Ramanujan polynomials, which 
have emerged in various contexts,  
see \cite{Berndt-I, Berndt, Chen-Guo-2001, Chen-Yang, Dumont-Ramamonjisoa, Guo, Ramanujan, Randazzo-2019, Shor, Wang-Zhou, Zeng, Guo-Zeng-2007, Lin-Zeng-2014, Sokal-2023, Chen-Sokal-2023, Josuat-Verges},
one should be reminded of the three intertwining scenes.
First, they refer to the polynomials
 $\psi_k(r,x)$ defined by Ramanujan 
 via the following relation \cite{Berndt-I, Ramanujan}.
Assume that $r$ is a nonnegative integer and that
\begin{equation}\label{equation:psi}
  \sum_{k=0}^\infty \frac{(x+k)^{r+k}e^{-u(x+k)}u^{k}}{k!}
  =\sum_{k=1}^{r+1}\frac{\psi_k(r,x)}{(1-u)^{k+r}},
\end{equation}
from which the following recurrence relation follows:
\begin{equation}\label{equation:relation_original}
  \psi_k(r+1,x)=(x-1)\psi_k(r,x-1)+\psi_{k-1}(r+1,x)-\psi_{k-1}(r+1,x-1),
\end{equation}
where  $\psi_1(0,x)=1$, $\psi_0(r,x)=0$ and
$ \psi_k(r,x)=0$ for $k>r+1$. 
The first few values
of $\psi_k(r,x)$ are given in Table \ref{Tab:1}.

\begin{table}[ht]
\centering
\caption{$\psi_k(r,x)$ for $0 \leq r \leq 3$ and 
$1 \leq k \leq r+1.$}
\label{Tab:1}
\vspace{3pt}
\begin{tabular}{|c|llll|c|}
\hline
\multicolumn{1}{|l|}{$r \backslash k$}& \multicolumn{1}{l}{$1$} & \multicolumn{1}{l}{$2$} & \multicolumn{1}{l}{$3$} & \multicolumn{1}{l}{$4$} &  \multicolumn{1}{|l|}{$\sum_k$} \\ 
\hline
$0$ & $1$ & $0$ & $0$ & $0$ &  $1$ \\
$1$ & $x-1$ & $1$ & $0$ & $0$ & $x$ \\
$2$ & $x^2-3x+2$ & $3x-5$ & $3$  & $0$ &  $x^2$ \\
$3$ & $x^3-6x^2+11x-6$ \quad \quad & $6x^2-26x+26$ \quad \quad & $15x-35$ \quad \quad & $15$ \quad &  $x^3$ \\ 
\hline
\end{tabular}
\end{table}

From the defining relation \eqref{equation:psi}, Berndt-Evans-Wilson 
\cite{Berndt-I,Berndt} derived
another recurrence relation
\begin{equation} \label{Berndt}
\psi_k(r,x)=(x-r-k+1)\psi_k(r-1,x)+(r+k-2)\psi_{k-1}(r-1,x).
\end{equation}

For the background of $\psi_k(r,x)$, we refer the reader to Ramanujan \cite{Ramanujan} and Berndt-Evans-Wilson \cite{Berndt}. 
Assume that $a$ is real with $|a|\ge e$. Then the 
solution of the equation 
\begin{equation} \label{xa}
    x=a\log x 
\end{equation} admits
the following expansion for any real number $n$,
\[
\frac{x^n}{n}=\sum_{k=0}^{\infty}\frac{(n +k)^{k-1}}{a^k k!}.
\] 
Ramanujan extended this sum to 
\[F_r(n)=\sum_{k=0}^{\infty}\frac{(n+k)^{r+k}}{a^k k!},\] 
and considered the question of how to determine the function $\psi_r(n)$ such that  
\begin{equation}  \label{psif}
    x^n \psi_r(n)=F_r(n). 
\end{equation}
When  $r$ is a nonnegative integer, by setting $x=e^u$, we see from (\ref{xa}) that $a=e^u/u$. Under this substitution, $\psi_r(n)$ is the left side of (\ref{equation:psi}). 
Ramanujan showed that for any nonnegative
integer $r$, 
\begin{equation}\label{equation:psi_sum}
\sum_{k=1}^{r+1}\psi_k(r,x)=x^r.
\end{equation}

 As remarked by Berndt-Evans-Wilson \cite{Berndt},
 no combinatorial perspective seems to be alluded to in Ramanujan’s original definition. Just on the contrary, the 
   polynomials $\psi_k(r,x)$ are
   destined to be under the
   umbrella of combinatorial analysis.
    
   On a different stage, Dumont-Ramamonjisoa \cite{Dumont-Ramamonjisoa}, with the aid of a 
   context-free grammar,  
    uncovered a combinatorial
   structure underneath an equation also considered by Ramanujan,
\begin{equation} \label{re}
x = y e^{-y} + \frac{a-1}{a} (e^{-y} -1). 
\end{equation}
They obtained expansions of
$y$ and $e^y$ as power series in $x$, where the coefficients are determined by the distribution of improper edges in labeled rooted trees. An associated insertion algorithm was also found.
Around the same time, Shor \cite{Shor} independently developed the same insertion algorithm.  Subsequently, Zeng \cite{Zeng} noted that these polynomials in Shor \cite{Shor} and Dumont-Ramamonjisoa \cite{Dumont-Ramamonjisoa} are linked to the Ramanujan polynomials $\psi_k(r,x)$, thereby providing a combinatorial interpretation.

Let us recall the notion of an improper edge of a labeled rooted tree. 
Let $n\geq 1$ and write $[n]=
\{1 ,2 , \ldots, n\}$. An edge of 
a rooted tree \( T \) on $[n]$
is represented by a pair \((i, j)\) of vertices with \( j \) being a child of \( i \). We say that  \((i, j)\)  is improper if there exists a descendant of \( j \) that is smaller than \( i \), bearing in mind that any vertex of \( T \) is considered as a descendant of itself; otherwise, \((i, j)\) is called a proper edge. We use
${\cal R}_{n,\,k}$ to denote the set 
of rooted trees on $[n]$ with $k$ improper edges, and
let 
$R(n,k)$  denote the number of rooted trees
in ${\cal R}_{n,\,k}$. Meanwhile, 
we use
${\cal T}_{n+1,\,k}$ to denote the set of rooted trees 
on $[n+1]$ with root $1$ and 
with $k$ improper edges, and let 
$T(n,k)$  denote
the cardinality of ${\cal T}_{n+1,\,k}$.  
The numbers $R(n,k)$ and $T(n,k)$  
are listed as sequences A054589 and A217922 in OEIS.

The following polynomials  are also referred to as the Ramanujan polynomials, 
\begin{align*}
R_n(u) & =\sum_{k= 0}^{n-1}  R(n,k) \, u^k, \\[6pt]
T_n(u) & = \sum_{k= 0}^{n-1} T(n,k)\,u^k .
\end{align*}
Their generating functions are defined by
\begin{align*}
    R(u,t) & = \sum_{n \geq 1} R_n(u) \frac{t^n}{n!},  
    \\[6pt]
    T(u,t) & =1+ \sum_{n \geq 1} T_{n}(u)\frac{t^n}{n!}.  
\end{align*}

As shown in \cite{Dumont-Ramamonjisoa}, (\ref{re}) gives rise to the 
following expansions,
\begin{align}
    y &  = \sum_{n \ge 1} \left( \sum_{k=0}^{n-1}
      R(n,k) a^{n+k} \right) \frac{x^n}{n!}, \label{yRnk}\\[6pt]
      e^y & =1+ \sum_{n \ge 1} \left(\sum_{k=0}^{n-1}
      T(n,k) a^{n+k} \right) \frac{x^n}{n!}. \label{eyTnk} 
\end{align}

Around the same time, Shor \cite{Shor} 
considered polynomials $Q_{n,\,k}(x)$ when
$x$ is treated as a positive
integer, which we also call the Ramanujan polynomials. 
Zeng \cite{Zeng} gave an interpretation of $Q_{n,\,k}(x)$, where 
$x$ is regarded as an indeterminate. The Ramanujan polynomials $Q_{n,\, k}(x)$ are given by
\begin{equation}\label{zeng-1} 
Q_{n,\, k}(x)=
\sum\limits_{T\in\mathcal{T}_{n+1,\,k}}x^{\deg_{T}(1)-1},
\end{equation}
or equivalently, 
\begin{equation} \label{zeng-2}
Q_{n,\,k}(x)
=\sum\limits_{T\in\mathcal{R}_{n,\,k}}(x+1)^{\deg_{T}(1)},
\end{equation}
where $\deg_T(i)$ denotes the degree of $i$ in $T$. The degree of a vertex $i$ in a
rooted tree $T$ is defined to be the number of children of $i$.
Zeng \cite{Zeng} also observed that 
\begin{equation} \label{Qnkx}
Q_{n,\,k}(x)=\psi_{k+1}(n-1,x+n). 
\end{equation}
Consequently,
Ramanujan's identity
\eqref{equation:psi_sum} aligns with Shor's identity \cite{Shor}: 
\begin{equation*}
\sum_{k=0}^{n-1}Q_{n,\,k}(x)=(x+n)^{n-1}.
\end{equation*}

In this context, the work of Dumont-Ramamonjisoa \cite{Dumont-Ramamonjisoa} can be viewed as an investigation into the special values of $Q_{n, \, k}(x)$. In particular, they found that $ Q_{n,\,k}(0)=R(n,k)$,
$Q_{n,\,k}(1)=T(n,k)$, and  $Q_{n,\,k}(-1)$ equals
the number of rooted trees on $[n]$ with $k$ improper edges for which the vertex $1$ is a leaf.
It is worth mentioning that   Wang-Zhou \cite{Wang-Zhou}
showed that 
$Q_{n,\,k}(-1)$ is related to refined orbifold Euler characteristics of the moduli space of stable curves of genus $0$ with $n$ marked points. 

This paper is organized as follows.
In Section 2, we present a  grammar for the
Ramanujan polynomials with an extra
label marking improper edges. We call this grammar literal because
each label
signifies exactly one of the cases in
the recursive construction of Shor and Dumont-Ramamonjisoa. 
Section 3 is devoted to a grammatical
calculus for the generating functions 
for the polynomials
$R_n(u)$, $T_n(u)$
and $Q_{n,\,k}(x)$, where constant properties
play a central role.
In Section 4, we 
provide a grammatical formulation of a bijection
in connection with the recurrence relation 
due to Berndt-Evans-Wilson and Shor. 

\section{A literal labeling scheme for rooted trees}

To provide a solution of the functional equation
(\ref{re}), Dumont-Ramamonjisoa   \cite{Dumont-Ramamonjisoa}
introduced the following grammar
\begin{equation}\label{D-D}
 A \to A^3S, \; S \to AS^2,
 \end{equation} and gave it a combinatorial interpretation. Let $D$ be the formal derivative associated with this grammar. For \( n \geq 1 \), Dumont-Ramamonjisoa obtained 
\begin{align*}
    D^{n-1}(AS) &  = R_n(A),
    \\[6pt]
     D^{n}(S)  & =  A^n S^{n+1}T_n(A).
\end{align*} 

In this section,
we present a labeling scheme for rooted trees (or planted rooted trees) based on the 
insertion algorithm due to 
  Shor \cite{Shor}  and 
Dumont-Ramamonjisoa \cite{Dumont-Ramamonjisoa}; see also
\cite{Randazzo-2019}. We call this labeling scheme
literal in the sense that 
each label corresponds to exactly
one of the cases in the recursive construction.
In principle, such an understanding of a labeling scheme may be instrumental 
in converting a recursive construction into a grammar. 
Moreover, the grammar generated this way might be
more convenient for carrying out the grammatical
calculus.

A rooted tree $T$ on $[n]$ can be regarded as 
a planted tree on $\{0,1,\ldots,n\}$ for which
$0$ is the root and the root has exactly one child.
Let $T$ be a rooted tree on $[n-1]$ with $n \geq 1$.
We describe  the
labeling scheme along with
the algorithm to insert the element $n$ into
$T$. For $n=0$, there is only one planted tree with the root $0$. 

\begin{itemize}
 
\item[1.] For any vertex $i$, $n$ may be added as a child of $i$. See Figure \ref{fig:z} for an illustration. 
           A new edge is formed, and a new vertex is formed. 
           If we use $z$ to label a vertex, and $v$ to label a new edge.
           Then we get a rule: 
           \[ z \to z(vz).\]
           The labels $v$ and
           $z$ in
           parentheses are meant to be
           the labels for the new vertex and the
           new edge (which is proper). The operation 
           in this case is called a $z$-insertion. 

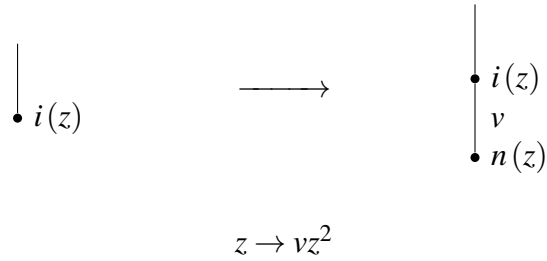
\begin{figure}[ht]
\centering
\begin{minipage}{0.3\textwidth}
\centering
\begin{tikzpicture}[scale=0.8,
  level 1/.style = { level distance   = 13mm,
                     sibling distance = 36mm },
  level 2/.style = { level distance   = 13mm,
                     sibling distance = 20mm },
  double line/.style={double, double distance=1pt, line width=0.18mm}
  ]
  
\node [ns,label=0:]{}[grow=down]
    child {node [nsf,label=0:{$i \, (z)$}](i){}
    };
   
\end{tikzpicture}
\end{minipage}%
\begin{minipage}{0.1\textwidth}
\centering
$\displaystyle \xrightarrow{\hspace*{1cm}}$ 
\end{minipage}%
\begin{minipage}{0.3\textwidth}
\centering
\begin{tikzpicture}[scale=0.8,
  level 1/.style = { level distance   = 13mm,
                     sibling distance = 36mm },
  level 2/.style = { level distance   = 13mm,
                     sibling distance = 36mm }, 
  level 3/.style = { level distance   = 13mm,
                     sibling distance = 10mm },
  double line/.style={double, double distance=1pt, line width=0.18mm}]
  
\node [ns,label=0:]{}[grow=down]
    child {node [nsf,label=0:{$i\,(z)$}](n){}
    child {node [nsf,label=0:{$n \,(z)$}](j){}}
    };
    
    \node [right=3mm, above=2.5mm] at (j){$v$};
\end{tikzpicture}
\end{minipage}

\vspace{5mm} 

\begin{minipage}{0.5\textwidth}
\centering
$z \to vz^2$
\end{minipage}

\vspace{5mm} 

\caption{A $z$-insertion.}
\label{fig:z}
\end{figure}

 \item[2.] For any edge (either proper or improper), 
            when $n$ is
           added, see Figure \ref{fig:v}, a new edge is created, which is improper. 
           So we should consider $v$ as a label for every edge. 
           Moreover, we should use an
           additional label $u$ to mark an improper edge. In this case, only edges labeled by $v$ are selected for the
           operation of inserting $n$. Notice that a new vertex is created.
           We get the rule: 
           \[ v \to v(uvz).\]
           The operation in this case is called a $v$-insertion, which does not change the degree of the existing vertices. 

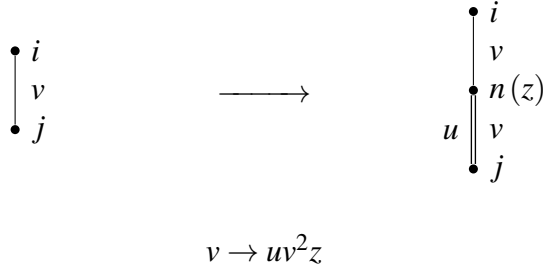
\begin{figure}[ht]
\centering
\begin{minipage}{0.3\textwidth}
\centering
\begin{tikzpicture}[scale=0.8,
  level 1/.style = { level distance   = 13mm,
                     sibling distance = 36mm },
  level 2/.style = { level distance   = 13mm,
                     sibling distance = 20mm },
  double line/.style={double, double distance=1pt, line width=0.18mm}
  ]
  
\node [nsf,label=0:$i$](i){}[grow=down]
    child {node [nsf,label=0:{$j$}](j){}
    };
    
    \node [right=3mm, above=2.5mm] at (j){$v$};
\end{tikzpicture}
\end{minipage}%
\begin{minipage}{0.1\textwidth}
\centering
$\displaystyle \xrightarrow{\hspace*{1cm}}$ 
\end{minipage}%
\begin{minipage}{0.3\textwidth}
\centering
\begin{tikzpicture}[scale=0.8,
  level 1/.style = { level distance   = 13mm,
                     sibling distance = 36mm },
  level 2/.style = { level distance   = 13mm,
                     sibling distance = 36mm }, 
  level 3/.style = { level distance   = 13mm,
                     sibling distance = 10mm },
  double line/.style={double, double distance=1pt, line width=0.18mm}]
  
\node [nsf,label=0:$i$](i){}[grow=down]
    child {node [nsf,label=0:{$n\,(z)$}](n){}
    child {node [nsf,label=0:{$j$}](j){}}
    };
    
    \node [right=3mm, above=2.5mm] at (n){$v$};
    \node [right=3mm, above=2.5mm] at (j){$v$};
    \node [left=3mm, above=2.5mm] at (j){$u$};
    \draw[double line] (n)--(j);
\end{tikzpicture}
\end{minipage}

\vspace{5mm} 

\begin{minipage}{0.5\textwidth}
\centering
$v \to uv^2z$
\end{minipage}

\vspace{5mm} 

\caption{A $v$-insertion.}
\label{fig:v}
\end{figure}

\item[3.] Each improper edge can be used to add $n$ in another way. See Figure \ref{fig:u} for an illustration. In this case, an improper
edge is created. Since a new vertex is
 created, we should use the rule: 
 \[ u \to u (uvz).\]
 The operation in this case is called a $u$-insertion and it does not change the degree of the existing vertices 
 except the vertex $i$. Figure 3 depicts the insertion, where we adopt the notation $\beta(j)$  for the minimum vertex among the vertices in the subtree rooted at $j$. 
\end{itemize}

\begin{figure}[ht]
\centering
\begin{minipage}{0.4\textwidth}
\centering
\begin{tikzpicture}[scale=0.8,
  level 1/.style = { level distance   = 13mm,
                     sibling distance = 36mm },
  level 2/.style = { level distance   = 13mm,
                     sibling distance = 20mm },
  double line/.style={double, double distance=1pt, line width=0.18mm},
  edge from parent/.style={ }
  ]
  
\node [ns,label=90:](root){}[grow=down]
    child {node [nsf,label=30:{$i$}](i){}
        child {node [nsf,label=270:{$j_1$}](1){}}
        child {node [nsf,label=270:{$j_d$}](d){}}
        child {node [nsf,label=270:{$j_{d+1}$}](s){}}
        child {node [nsf,label=270:{$j_l$}](l){}}
    };
    \coordinate (left) at ([xshift=-0.55cm,yshift=0.18cm] s);
    \coordinate (right) at ([xshift=0.5cm,yshift=-0.15cm] l);
    \draw (left) rectangle (right);
    \node [right=4mm] at (1){$\ldots$};
    \node [right=3.5mm] at (s){$\ldots$};
    \node [left=-1mm, above=2.5mm] at (d){$u$};
    \node [right=24mm, below=8mm] at (1){$\beta(j_1) < \beta(j_2) < \cdots < \beta(j_l)$};
    \draw[double line] (i) -- (1);
    \draw[double line] (i) -- (s);
    \draw[double line] (i) -- (d);
    \draw (root) -- (i);
    \draw (i) -- (l);
\end{tikzpicture}
\end{minipage}%
\begin{minipage}{0.1\textwidth}
\centering
$\displaystyle \xrightarrow{\hspace*{1cm}}$ 
\end{minipage}%
\begin{minipage}{0.4\textwidth}
\centering
\begin{tikzpicture}[scale=0.8,
  level 1/.style = { level distance   = 13mm,
                     sibling distance = 36mm },
  level 2/.style = { level distance   = 13mm,
                     sibling distance = 20mm }, 
  level 3/.style = { level distance   = 13mm,
                     sibling distance = 10mm },
  double line/.style={double, double distance=1pt, line width=0.18mm}]
  
\node [ns,label=90:]{}[grow=down]
    child {node [nsf,label=0:{$n \, (z)$}](n){}
        child {node [nsf,label=270:{$j_1$}](1){}}
        child {node [nsf,label=270:{$j_d$}](d){}}
        child {node [nsf,label=0:{$i$}](i){}
        child {node [nsf,label=270:{$j_{d+1}$}](s){}}
        child {node [nsf,label=270:{$j_{l}$}](l){}}
        }
    };
    
    \node [right=5mm] at (1){$\ldots$};
    \node [right=0.5mm] at (s){$\ldots$};
    \node [left=9mm, above=1mm] at (i){$u$};
    \node [right=-3mm, above=2.5mm] at (i){$v$};
    \node [left=2mm, above=2.5mm] at (d){$u$};
    \draw[double line] (n) -- (1);
    \draw[double line] (n) -- (i);
    \draw[double line] (n) -- (d);
    \draw[double line] (i) -- (s);
    \coordinate (left) at ([xshift=-0.4cm,yshift=0.15cm] s);
    \coordinate (right) at ([xshift=0.4cm,yshift=-0.15cm] l);
    \draw (left) rectangle (right);
\end{tikzpicture}
\end{minipage}

\vspace{5mm} 

\begin{minipage}{0.5\textwidth}
\centering
$u \to u^2vz$
\end{minipage}

\vspace{5mm} 

\caption{A $u$-insertion.}
\label{fig:u}
\end{figure}
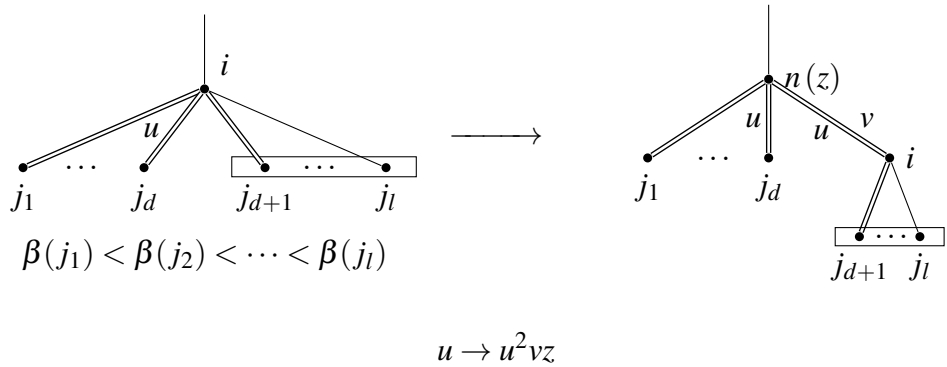

It can be readily seen that
the insertion algorithm is reversible. 
In accordance with the above procedure, it is natural
to introduce the following labeling scheme for a planted tree $T$ on $\{0,1,\ldots,n\}$ for any $n\geq 1$:
\begin{itemize}
         
    \item each vertex $i$ other than $0$ is labeled by $z$;

    \item each edge is labeled by $v$;

    \item each improper edge carries an additional label $u$. 
\end{itemize}
In other words, a proper edge is labeled by $v$, whereas an improper edge is labeled by $uv$.  

It should be noted that
this labeling scheme has a special
feature that the labels are in one-to-one correspondence with the insertions of the
recursive construction. 
 In this setting, an extra label 
is at our disposal. This updated labeling scheme  
seems to offer more 
flexibility in further exploring the Ramanujan polynomials
subject to certain constraints.

The rules are summarized into the following grammar:
\begin{equation} \label{R-G} 
G=\{z\to vz^2, \; v\to uv^2z, \; 
     u\to u^2vz \}. 
\end{equation}
Let $D$ denote the formal
derivative with respect to this grammar $G$.

\begin{theorem}  
For $n\geq 1$, we have
\[ D^{n-1}(vz) = v^{n}z^n R_n(u).\]
\end{theorem}

Likewise, the polynomials 
$T_n(u)$ can be generated by applying the
operator $D^{n-1}$ to $z$.

\begin{theorem}
\label{thm: D-T}
For $n\geq 1$, we have
\[ D^{n-1}(z) = v^{n-1} z^{n} T_{n-1}(u) , \]
where $T_0(u)$ is defined to be $1$.
\end{theorem}

For $0 \leq n \leq 4$,
the  values of $D^n(vz)$ and $D^n(z)$ are given below,
\begin{align*}
    D^0(vz) & =vzR_1(u) = vz, \\[3pt]
    D^1(vz) & =v^2z^2R_2(u) = v^2z^2(1+u),\\[3pt]
    D^2(vz) & = v^3z^3R_3(u) =v^3z^3(2+4\, u +3\, u^2), \\[3pt]
    D^3(vz) & =v^4z^4R_4(u) =  v^4z^4(6+18\,  u+25\, u^2 +15\,  u^3),\\[3pt]
    D^4(vz) & =v^5z^5R_5(u) =  v^5z^5(24+96\,  u+190\, u^2 +210\, u^3+105\, u^4)
\end{align*}
and
\begin{align*}
    D^0(z)  & = zT_0(u)=z, \\[3pt]
    D^1(z) &  =vz^2T_1(u)= vz^2, \\[3pt]
  D^2(z) & =v^2z^3T_2(u)= v^2z^3(2 + u), \\[3pt]
  D^3(z) & =v^3z^4T_3(u)=v^3z^4(6 +7 u +3 u^2), \\[3pt]
   D^4(z)& =v^4z^5T_4(u)=v^4z^5(24+46\,u +40\,u^2+15\,u^3). 
\end{align*}

\section{The grammatical calculus}

In this section, we develop a grammatical
calculus leading to functional equations for the generating 
functions of $R_n(u)$,  $T_n(u)$ and $Q_{n,\,k}(x)$.

Recall that for a Laurent series 
$f$ in the variables involved in the grammar $G$, the associated generating function is defined by
\[ \gen (f, t) = \sum_{n=0}^{\infty} D^n(f) \frac{t^n}{n!}. \]
For two Laurent series $f$ and $g$, the Leibniz formula holds for $n\geq 0$,
\[ D^n(fg) = \sum_{k=0}^n \binom{n}{k}
D^k(f) D^{n-k}(g),\]
which implies the multiplicative property 
\[ \gen (fg, t) = \gen (f , t)\,  \gen ( g, t).   \]
In particular, we have
\[ \gen(f,t) \gen(f^{-1},t) =1. \]

Constants and eigenfunctions are 
basic ingredients of the grammatical calculus,
see, for example,  \cite{Ji, Fu, 
Chen-Fu-2023-A, Chen-Fu-2023-B}. A function $f$
is called a constant with respect to $D$ if $D(f)=0$, and a $k$-th order constant if $D^k(f)=0$. 
An eigenfunction $f$ is a function satisfying $D(f)=cf$, where $c$ is a constant.

First, we consider the 
functional equation satisfied by the
generating function of $T_n(u)$. 
  
\begin{theorem}\label{Gen-z} Given the grammar 
$$
G=\{z\to vz^2, \; v\to uv^2z, \; 
     u\to u^2vz \},
$$
the generating function $\gen(z,t)$ satisfies 
\begin{align}\label{gen-z}
     \gen(z,t)=ze^{(z^{-1}-u^{-1}z^{-1}+u^{-1}vt)\gen(z,\,t)+u^{-1}-1}.
\end{align}
\end{theorem}

From Theorem \ref{thm: D-T}, we have for $n \ge 0$, 
\begin{align*}
    D^{n}(z) \mid_{v=z=1}=T_{n}(u).
\end{align*}
Thus we obtain the following functional equations for $T(u,t)$ by setting  \(v=1\) and \(z=1\) in Theorem \ref{Gen-z}: 
\begin{equation}
    \label{eq-tt}
     T(u,t)= e^{(1-u^{-1} +u^{-1} t)T(u,\,t)+u^{-1}-1}.
\end{equation}

For a labeled tree $T$ on $[n+1]$ rooted at 
$1$, any edge incident to the root is proper. 
After relabeling, such a rooted tree can be viewed as 
a forest $F$ on $\{1,2,\ldots, n\}$, and an edge of 
$T$ is improper if and only if it is improper in a component of $F$.
It follows that 
\[ T(u,t) = e^{R(u,\,t)}. \]
Combining this identity with  (\ref{eq-tt}), we obtain
\[ e^{R(u,\,t)} = 
       e^{(1-u^{-1} +u^{-1} t) e^{R(u,\,t)}+u^{-1}-1},
\]
and hence 
\begin{equation}  \label{eq-rt}
R(u,t)   = 
        {(1-u^{-1} +u^{-1} t)e^{R(u,\,t)}+u^{-1}-1}. 
\end{equation}

Equation \eqref{eq-rt} can be transformed into
 (\ref{re}) of Ramanujan under the substitutions $u=a$, $t=ax$ and $y=R(a,ax)$. 

The following properties  are needed  in the 
proof of Theorem \ref{Gen-z}. Since the formal derivative $D$ with respect to grammar $G$ can be viewed as
a differential operator
\[ D= u^2 v z \frac{\partial}{ \partial u} + 
      uv^2z \frac{\partial}{ \partial v}
     + vz^2 \frac{\partial}{ \partial z},\]
and so these properties can be derived by
solving the corresponding partial differential equations.

\begin{proposition}
\label{prop:3.2}
The following relations hold
    \begin{align}
    D(u^{-1}v)&=0,\label{constant-1}\\[3pt]
    D(ze^{u^{-1}})&=0,\label{z1}\\[3pt]
    D((u-1)v^{-1}z^{-1})&=1, \label{z1a} \\[3pt]
    D(z^{-1}-u^{-1}z^{-1})&=u^{-1}v, \label{z2} \\[3pt]
    D(e^{u^{-1}}(uv^{-1}-v^{-1}))&=ze^{u^{-1}}. \label{z3}
\end{align}
\end{proposition}

These identities can be verified by direct computations. Especially, we note that
\eqref{z2}  follows from \eqref{constant-1}  and \eqref{z1a} since
\begin{align*}
    z^{-1}-u^{-1}z^{-1} &=(u^{-1}v)\cdot \left((u-1)v^{-1}z^{-1}\right),
\end{align*}
and similarly, \eqref{z3} follows from \eqref{z1} and \eqref{z1a},
since 
\begin{align*}
    e^{u^{-1}}(uv^{-1}-v^{-1})=(ze^{u^{-1}})\cdot \left((u-1)v^{-1}z^{-1}\right).
\end{align*}

We are now ready to give a grammatical
derivation of Theorem \ref{Gen-z}. 

\noindent
{\it Proof of Theorem \ref{Gen-z}.} 
In view of
\eqref{z1}, we see that
\begin{align}
    \label{eq1}
    \gen(ze^{u^{-1}},t) 
     =ze^{u^{-1}}.
\end{align}
Noting  that  
\begin{align*}
    \gen(e^{u^{-1}},t)=e^{\gen(u^{-1},\,t)},
\end{align*}
by (\ref{eq1}), we obtain 
\begin{equation}
 \gen(z,t)e^{\gen(u^{-1},\,t)} = ze^{u^{-1}},
\end{equation}    
and thus
\begin{align}
    \begin{aligned}\label{general-uz}
    e^{\gen(u^{-1},\,t)}=ze^{u^{-1}}\gen(z^{-1},t).
    \end{aligned}
\end{align}
By \eqref{z2} and \eqref{constant-1}, we find that
\begin{equation}
    \label{eq2}
   \gen(z^{-1}-u^{-1}z^{-1},t) 
   = z^{-1}-u^{-1}z^{-1} + 
   u^{-1}v t. 
\end{equation}
Therefore,
\begin{equation}
     \gen(z^{-1},t)-\gen(u^{-1},t)\gen(z^{-1},t) 
     =z^{-1}-u^{-1}z^{-1}+u^{-1}vt.
\end{equation}
Multiplying both sides by $\gen(z,t)$, we get
\begin{equation}
\label{eq:3.14}
\gen(u^{-1},t)=\gen(z,t)(\gen(z^{-1},t)- z^{-1}+u^{-1}z^{-1}-u^{-1}vt).
\end{equation}
Substituting \eqref{eq:3.14} into \eqref{general-uz}, we obtain
\begin{align}
    e^{\gen(z^{-1},\,t)}=e^{z^{-1}-u^{-1}z^{-1}+u^{-1}vt}(ze^{u^{-1}}\gen(z^{-1},t))^{\gen(z^{-1},\,t)},
\end{align}
which implies (\ref{gen-z}). This completes 
the proof. \qed 

The functional equation 
for $\gen(u,t)$ can
be established by using
 \eqref{gen-z} and \eqref{general-uz}. 
More precisely, 
we have
\begin{align}\label{gen-u}
    (1-\gen(u^{-1},t))e^{\gen(u^{-1},\,t)}=e^{u^{-1}}(1-u^{-1} +u^{-1}vzt).
\end{align}
As for $\gen(v,t)$, by \eqref{constant-1} we have
\begin{align*}
    \gen(u^{-1}v,\,t)=\gen(u^{-1},t)\gen(v,t)=u^{-1}v.
\end{align*}
Thus,
\begin{align*}
    \gen(u^{-1},t)=u^{-1}v\,\gen(v^{-1},t),
\end{align*}
which, together with  \eqref{gen-u}, yields
\begin{align}
    \left(1-u^{-1}v\,\gen(v^{-1},t)\right)e^{u^{-1}v\,\gen(v^{-1},t)}=e^{u^{-1}}(1-u^{-1} +u^{-1}vzt).
\end{align}

We next turn to 
the grammatical derivation of
the functional equation of the generating function of $Q_{n,\,k}(x)$. 
In doing so, we introduce two additional labels $a$ and $x$ along with their corresponding substitution rules, on top of  the grammar in Section 2:
$$
G=\{z\to vz^2, \; v\to uv^2z, \; 
     u\to u^2vz \}.
$$

Given a rooted tree $T$ on $[n]$ with root $1$, label the root by $a$. The children of the root are labeled by $x$ and the rest of the vertices are labeled by $z$. The edges are labeled by $v$, and the improper edges have an extra label $u$.

Starting with the root $1$ with label $a$, $D(a)$ gives a 
labeled tree on $[2]$ with $2$ being a child of the root. 
This gives the rule: 
\[ a \to  axv.\] 
Inspecting the insertions 
involving the children of the
root and other vertices, we get the rules
\[ x \to xv  z , \;\, z  \to v z ^2,\; v \to uv^2z,\; u  \to u^2v z, \]
and so we have the grammar
\begin{align}
\label{G-Q-our}
G=\{a \to  axv,\; x\to xvz, \;z\to  vz^2, \; v\to uv^2z, \; 
     u\to u^2vz \}. 
\end{align}
By a slight abuse of notation, we still denote the updated grammar by $G$.
This grammar has an extra label compared with the grammar in \cite{Chen-Yang}.

\begin{theorem}
\label{thm:gen-a}
For $n \ge 1$, we have
\begin{align}
    D^n(a)=av^n\sum_{k=0}^{n-1}\;\sum_{T \in \mathcal{T}_{n+1,\,k}}x^{\deg_T(1)}z^{n-\deg_T(1)}u^k,
\end{align}
or equivalently,
\begin{align}
    D^n(a)&=axv^nz^{n-1}\sum_{k=0}^{n-1} Q_{n,\,k}(xz^{-1})u^k. \label{Dna}
\end{align}
\end{theorem}

In the spirit of the grammatical calculus for the generating functions of $R_n(u)$ and $T_n(u)$, we can carry out a grammatical
calculus for the generating function of 
$Q_{n,\,k}(x)$, leading to the functional
equation \eqref{re} of Ramanujan. 

First, we find it more convenient to
make a change of variables by setting $b=xv$ and $c=vz$
for the grammar $G$ in \eqref{G-Q-our}. Let $D_G$ 
denote the formal derivative with respect to $G$.
Since 
\begin{align}
    D_G(xv)=xv^2z(1+u)=bc(1+u),\\[3pt]
    D_G(vz)=v^2z^2(1+u)=c^2(1+u),
\end{align}
the grammar $G$ is transformed into 
\begin{align}
\label{H}
    H=\{a \to ab,\, b \to bc(1+u),\, c \to c^2(1+u),\, u \to cu^2\}.
\end{align}
Let $D_H$ denote the formal derivative with
respect to $H$. Evidently, 
\begin{align}
\label{H-G}
D_H(a) \mid _{b=xv,\,c=vz} \, =\,D_G(a).
\end{align}

In the notation of the grammar 
$H$, the following properties hold.

\begin{proposition}
\label{prop:a}
We have 
    \begin{align}
        D_H(uc^{-1}e^{-u^{-1}})&=0, \label{ac-1}\\
        D_H(bc^{-1})&=0, \label{ac-2}\\
        D_H(ae^{bc^{-1}u^{-1}})&=0, \label{ac-3}\\
        D_H((u-1)c^{-1})&=1. \label{ac-4}    
    \end{align}   
\end{proposition}

The following theorem gives a functional
equation satisfied by the generating function 
$\gen(a,t)$. 

\begin{theorem}\label{Zeng-equiv}
Let $y \in \mathbb{C}[[t,u]]$ be the solution of the equation 
\begin{align}
\label{eq:3.30}
    (1-u^{-1}+y)ue^{-y}=u-1+vzt.
\end{align}
Then
\begin{align}
   \gen(a,t)=ae^{xz^{-1}y}.\label{5.13}
\end{align} 
\end{theorem}

 \proof
Let $\gen_H$ denote the generating function with respect to grammar $H$ in \eqref{H}. 
It follows from \eqref{ac-2} and \eqref{ac-3} that
\begin{align}
    \gen_H(bc^{-1},t)&=bc^{-1},\label{ac-5}\\[6pt]
    \gen_H(a,t)e^{\gen_H(bc^{-1},\,t)
    \gen_H(u^{-1},\,t)}&=ae^{bc^{-1}u^{-1}}. \label{ac-6}
\end{align}
Substituting \eqref{ac-5} into \eqref{ac-6}, we get 
\begin{align}
\gen_H(a,t)=ae^{bc^{-1}(u^{-1}-\gen_H(u^{-1},\,t))} \label{ac-7}.
\end{align}
To compute $\gen_H(u^{-1},t)$, 
using the constant property
\eqref{ac-1}, we find that
\begin{align} \gen_H(u,t)\gen_H(c^{-1},t)e^{\gen_H(-u^{-1},\,t)}=uc^{-1}e^{-u^{-1}}, \label{ac-8}
\end{align}
and thus
\begin{align}
\label{eq:3.35}
\gen_H(c^{-1},t)=uc^{-1}e^{-u^{-1}}\gen_H(u^{-1},t)e^{\gen_H(u^{-1},\,t)}.
\end{align}
From \eqref{ac-4}, we obtain
\begin{align} 
    \gen_H(u-1,t)\gen_H(c^{-1},t)=(u-1)c^{-1}+t. \label{ac-9}
\end{align}   
Substituting \eqref{eq:3.35} into \eqref{ac-9} yields
\begin{align}
    ue^{-u^{-1}}e^{\gen_H(u^{-1},\,t)}(1-\gen_H(u^{-1},t))=u-1+ct.
\end{align}
Writing $\gen_H(u^{-1},t)=u^{-1}-y$, we have
\begin{align}
    ue^{-u^{-1}}e^{u^{-1}-y}(1-u^{-1}+y)=u-1+ct,
\end{align}
which simplifies to  
\begin{align}\label{ac-y}
    (1-u^{-1}+y)ue^{-y}=u-1+ct.
\end{align}
Hence by \eqref{ac-7}, we get
\begin{align}
    \gen_H(a,t)=ae^{bc^{-1}y},
\end{align}
where $y$ is the solution of \eqref{ac-y}.

Back to the grammar $G$, we conclude 
that  
\begin{align}
    \gen(a,t)=ae^{xz^{-1}y},
\end{align}
where $y \in \mathbb{C}[[t,u]]$ is the solution of the equation 
\begin{align}
    (1-u^{-1}+y)ue^{-y}=u-1+vzt.
\end{align}
This completes the proof. 
 \qed

In light of
Theorem \ref{thm:gen-a}, we see that 
\begin{align}
\label{eq:Gen-a}
    \gen(a,t)=a+ axz^{-1}\sum_{n\ge 1} \sum_{k=0}^{n-1} Q_{n,\,k}(xz^{-1})u^k \frac{(vzt)^n}{n!}.
\end{align}  
Note that Zeng  \cite{Zeng} defined the generating function of $Q_{n,\,k}(x)$ in the following form
    \begin{align}
        Y(u,t)=\sum_{n \ge 1}\sum_{k=0}^{n-1}\frac{Q_{n,\,k}(x)}{(1-u)^k}\frac{t^n}{n!}
    \end{align}
and showed that
    \begin{align} \label{yut}
        Y(u,t)=\frac{e^{xy}-1}{x},
    \end{align}
    where $y \in \mathbb{C}[[t,u]]$ is the solution of the equation 
\begin{align}\label{yut-2}
   (1-u)t=ye^{-y}+u(e^{-y}-1).
\end{align}

Notice that by setting $a=v=z=1$ and replacing
$u$ with  $1/(1-u)$, the functional equation of $Y(u,t)$ is equivalent to that of $\gen(a,t)$.

To conclude this section, 
we notice that equation \eqref{yut-2} is equivalent to Ramanujan's equation \eqref{re} 
by replacing $(1-u)t$ with $x$ and $u$ with $(a-1)/a$.
In fact, 
from the  expansion of $y$
given by Dumont-Ramamonjisoa   \eqref{yRnk},
resorting to the exponential formula for labeled structures, see Stanley \cite{Stanley}, 
we are led to the following expansion 
\begin{align} \label{exy}
        e^{xy} =1+x\sum_{n \ge 1}\left(\sum_{k=0}^{n-1} Q_{n,\,k}(x)  a^{n+k}\right) \frac{x^n}{n!},
\end{align} 
where $Q_{n,\,k}(x)$ are endowed with the
combinatorial interpretation as in \eqref{zeng-1}.

\section{The grammatical calculus behind a bijection}

This section is concerned with a grammatical
formulation of a recurrence
relation of the Ramanujan polynomials. 
Recall that
Berndt-Evans-Wilson \cite{Berndt} derived the 
recurrence relation \eqref{Berndt} for $\psi_k(r,x)$. On the other hand,
Shor  \cite{Shor} asked for a combinatorial 
proof of a recurrence relation for 
$Q_{n,\,k}(x)$, that is, 
\begin{equation}
\label{Berndt-Q}
Q_{n,\,k}(x) = (x - k + 1)Q_{n-1,\,k}(x+1) + (n + k - 2)Q_{n-1,\,k-1}(x+1).
\end{equation}
It turns out that these two 
recurrences are equivalent. Chen-Guo \cite{Chen-Guo-2001}
presented a rather involved bijection for this recurrence,
and later Guo \cite{Guo} found a simpler construction. Based on the grammar mentioned in Section 3, Chen-Yang \cite{Chen-Yang}
gave a grammatical proof using generating functions. Albeit these efforts,
a better combinatorial understanding is still in demand. 
Perhaps an alternative
insertion algorithm is needed for 
this purpose. 

Chen-Guo \cite{Chen-Guo-2001} showed that the recurrence (\ref{Berndt-Q}) 
can be deduced from the following correspondence. 

\begin{theorem} \label{Theorembij} For $n\geq 2$,
there is a bijection 
    \[ \mathcal{T}_{n,\,k}[\deg_T(2)>0,\deg_T(1)=r] \longleftrightarrow 
\mathcal{T}_{n,\,k+1}[\deg_T(n) >0,\deg_T(1)=r],\]
where the conditions in brackets specify degree constraints. 
\end{theorem} 

Let's translate the above
correspondence into a grammatical statement, 
where the grammar $G$ is given in  (\ref{G-Q-our}), that is,
\begin{align*}
G=\{a\rightarrow axv, \; x\rightarrow xvz, \; z\rightarrow vz^2, \; v\rightarrow uv^2z, \; 
     u\rightarrow u^2vz \}.
\end{align*} 
The insertion process of generating rooted trees on $[n]$ with root $1$ yields
\[D^{n-1}(a)=av^{n-1}\sum_{k=0}^{n-2}\;\sum_{T \in \mathcal{T}_{n,\,k}}x^{\deg_T(1)}z^{n-\deg_T(1)-1}u^k.\]
Start with the edge $(1,2)$, where the vertex 
$1$ is labeled by $a$, the vertex $2$ is labeled
by $x$ and the edge is labeled by $v$. Observe that once a vertex becomes an internal vertex, it remains an internal vertex. The condition $\deg(2)>0$ indicates that the vertex $2$ is not allowed to be a leaf. 
Concerning the trees on $[n]$ with $\deg(2)=0$, since the vertex $2$ never gets involved in an $x$-insertion, the generating polynomial of such trees is $xD^{n-2}(av)$. 

On the other hand, the condition $\deg(n)=0$ says that the vertex $n$ is a leaf.  Assume that $T$ is obtained from $T'$ with $n-1$ vertices by the insertion of a leaf $n$. There are three cases. If $n$ is a child of the root $1$, that is, an $a$-insertion gets involved, the generating polynomial will be $xvD^{n-2}(a)$. For an $x$-insertion or a $z$-insertion, the generating polynomials of both cases shall be
$(n-2)vzD^{n-2}(a)$. 

To give a grammatical restatement of the bijection in Theorem \ref{Theorembij}  the
number of 
 improper edges should be taken into consideration. This is reflected
 by  a factor $u$ on the left side
 in the following identity,  
\[u(D^{n-1}(a)-xD^{n-2}(av))
=D^{n-1}(a)-xvD^{n-2}(a)-(n-2)vzD^{n-2}(a).\]
So the grammatical identity we wish to establish
can be formulated as follows, which is not
hard to justify by induction.

\begin{theorem}
For $n \ge 2$, we have
\begin{equation}
    \label{grammar-recur}
    D^{n-2}(av) = \frac{u-1}{xu}D^{n-1}(a) + \frac{v}{u}\left(1 + (n-2)\frac{z}{x}\right)D^{n-2}(a).
\end{equation}
\end{theorem}

\proof
We proceed by induction on $n$. First, 
we observe the following constant properties 
\[
D\left(\frac{v}{u}\right) = 0, \quad D\left(\frac{z}{x}\right) = 0,\]
and 
\[D\left(\frac{u-1}{xu}\right) = \frac{vz}{xu}.
\]

Applying the operator $D$ to both sides of \eqref{grammar-recur}, by the Leibniz rule we find that
\begin{align*}
D^{n-1}(av) &= D\left( \frac{u-1}{xu}D^{n-1}(a) \right) + D\left( \frac{v}{u}\left(1 + (n-2)\frac{z}{x}\right)D^{n-2}(a) \right) \\[5pt]
&= \frac{vz}{xu}D^{n-1}(a) + \frac{u-1}{xu}D^{n}(a) + \frac{v}{u}\left(1 + (n-2)\frac{z}{x}\right)D^{n-1}(a).
\end{align*}
Thus,
\begin{align*}
D^{n-1}(av)= \frac{u-1}{xu}D^{n}(a) + \frac{v}{u}\left(1 + (n-1)\frac{z}{x}\right)D^{n-1}(a),
\end{align*}
and so the proof is complete by induction.
\qed

\vskip 6mm \noindent{\large\bf Acknowledgment.} 
 We wish to thank the referees for their insightful 
 comments.

\end{document}